\documentclass[11pt]{article}

\usepackage{amsfonts}
\usepackage{amssymb}
\usepackage{amsmath}
\usepackage{fullpage}
\usepackage{color}
\usepackage{slashbox}
\usepackage{amsthm}
\usepackage{dsfont}
\usepackage{graphicx}

\def\R{\mathbb{R}}

 \numberwithin{equation}{section}

\newtheorem{theorem}{Theorem}[section]

\newtheorem{remark}[theorem]{Remark}

\DeclareMathOperator*{\argmin}{\arg\!\min}

\begin{document}
\title{Optimal Point Sets for Total Degree Polynomial Interpolation in Moderate Dimensions}
\author{Max Gunzburger$^1$ and Aretha L. Teckentrup$^1$}
\date{}
\maketitle

\begin{center}
\begin{footnotesize}

\vspace{-0.75cm}

\noindent
${}^1$ Dept of Scientific Computing, Florida State University, 400 Dirac Science Library, Tallahassee FL 32306-4120 \\\indent {\tt gunzburg@fsu.edu}, {\tt ateckentrup@fsu.edu}

$  $
\end{footnotesize}
\end{center}

\abstract{This paper is concerned with Lagrange interpolation by total degree polynomials in moderate dimensions. In particular, we are interested in characterising the optimal choice of points for the interpolation problem, where we define the optimal interpolation points as those which minimise the Lebesgue constant. We give a novel algorithm for numerically computing the location of the optimal points, which is independent of the shape of the domain and does not require computations with Vandermonde matrices. We perform a numerical study of the growth of the minimal Lebesgue constant with respect to the degree of the polynomials and the dimension, and report the lowest values known as yet of the Lebesgue constant in the unit cube and the unit ball in up to 10 dimensions.
}

\section{Introduction}
One of the fundamental problems in approximation theory is the Lagrange interpolation problem: given a continuous function $f$, a set of interpolation points $X$ and a polynomial interpolation space $\Pi$, find an element of $\Pi$ which is equal to $f$ at the interpolation points $X$. Although the statement of this problem is very simple, many open questions remain on this topic, especially in the case of multivariate interpolation. For an overview of some of the current research topics in the area of multivariate interpolation, such as ensuring well-posedness of the Lagrange interpolation problem, the explicit construction of the interpolating polynomial for certain interpolation points and the derivation of remainder formulas, we refer the reader to the works \cite{gs00a,gs00b,cg10} and the references therein.

The question that will be addressed in this paper, is how to choose the interpolation points $X$ such that the interpolation error in an arbitrary continuous function $f$ is small. More precisely, we are looking for the set of interpolation points which gives the smallest possible upper bound on the interpolation error in an arbitrary continuous function. Classically, this optimal choice of interpolation points is given by those which minimise the Lebesgue constant (see section \ref{sec:prob} for more details). Minimising the Lebesgue constant amounts to solving a large scale non-linear optimisation problem, to which the true solution is not explicitly known, even in the case of univariate interpolation. 

With a particular focus on the multivariate setting, the aim of this paper is to study the structure of the optimal point sets resulting from a minimisation of the Lebesgue constant and also the behaviour of the corresponding Lebesgue constants themselves. Although there have been several studies on this problem in the bivariate setting (see for example \cite{vhs14,bsv12} and the references therein), not much is known in the general multivariate case. By employing a novel algorithm to numerically minimise the Lebesgue constant, we are able compute optimal point sets for interpolation in up to 10 variables. Our algorithm does not require the computation of Vandermonde determinants or the solution of linear systems of equations. The Lebesgue constants we compute are the lowest known values as yet, also in the bivariate setting already considered in \cite{vhs14,bsv12}. We further provide suggested growth rates of the Lebesgue constant which concur with current theoretical results and conjectures.

The structure of the remainder of this paper is as follows: In section \ref{sec:prob}, we give a mathematical formulation of the problem of finding optimal interpolation points in the general multivariate setting, and review relevant results from the literature. We give an algorithm that efficiently computes the optimal interpolation points in section \ref{sec:alg}, before we in section \ref{sec:num} present our findings on optimal interpolations points in various geometries in up to 10 dimensions. Some conclusions and suggestions for further work are given in section \ref{sec:conc}.

\section{Problem formulation}\label{sec:prob}
Let $\Pi_n^d$ be the space of all polynomials in $d$ variables with total degree less than or equal to $n$. The dimension of this space is $N = \binom{n+d}{d}$. Given a compact subset $D \subset \R^d$, and a set of distinct points $X_{n,d} = \{\xi_j\}_{j=1}^N$ in $D$, the Lagrange interpolation problem is then the following: for each $f \in C^0(D)$, find a polynomial $p_f \in \Pi_n^d$, such that
\begin{equation}\label{eq:int_con}
p_f(\xi_j) = f(\xi_j), \qquad j=1,\dots,N.
\end{equation}
If a unique such polynomial exists for all $f \in C^0(D)$, the point set $X_{n,d}$ is said to be {\em unisolvent}, and the polynomial $p_f$ is called the {\em Lagrange interpolant} of the function $f$. For the remainder of this paper, we will assume that the point set $X_{n,d}$ is unisolvent. In the case of univariate interpolation, the unisolvence of the set $X_{n,1}$ is equivalent to the points $\{\xi_j\}_{j=1}^N$ being distinct. In higher dimensions, the issue of unisolvence is more complicated and depends on the geometric configuration of the interpolation points. In general the set $X_{n,d}$ is unisolvent if, and only if, its points do not lie on a hyper surface of degree $n$, or equivalently, the Vandermonde determinant formed by the points is non--zero (cf \eqref{eq:lebpol_d}).

Let now $f$ denote an arbitrary function in $C^0(D)$. In general, the Lagrange interpolant of $f$ can be constructed by choosing a basis for $\Pi_n^d$, and subsequently determining the coefficients of $f$ in this basis using the interpolating conditions \eqref{eq:int_con}. A particular basis often used for this purpose are the {\em Lagrange fundamental polynomials}, denoted by $\{l_j\}_{j=1}^N $, which have the property
\[
l_j(\xi_k) = \delta_{jk}, \qquad j,k=1, \dots, N,
\]
where $\delta_{jk}$ denotes the Kronecker delta. The Lagrange interpolant of $f$ at any point $x \in D$ is then given by
\[
L_{n} f (x) = \sum_{j=1}^N f(\xi_j) l_j(x).
\]
An explicit expression for the Lagrange fundamental polynomials is given in terms of Vandermonde determinants. For any set of points $Y = \{y_j\}_{j=1}^N \subset D$, we denote by $V(Y)$ the $N \times N$ generalised Vandermonde matrix with entries $(V(Y))_{ij} = \phi_j(y_i)$, where $\{ \phi_j\}_{j=1}^N$ is any basis of $\Pi_n^d$, such as the monomial basis. The Lagrange fundamental polynomials are then given by 
\begin{equation}\label{eq:lebpol_d}
l_j(x) = \frac{\textrm{det} [V(\xi_1, \dots, \xi_{j-1}, x, \xi_{j+1}, \dots, \xi_N)]}{\textrm{det} [V(\xi_1, \dots, \xi_N)]}, \qquad j=1, \dots, N.
\end{equation} 
It is clear from \eqref{eq:lebpol_d} that the Lagrange fundamental polynomials can be constructed only if the Vandermonde determinant $\textrm{det} [V(\xi_1, \dots, \xi_N)$ is non--zero, which is equivalent to the points $X_{n,d}$ being unisolvent. 
In the special case of univariate interpolation, the expression \eqref{eq:lebpol_d} reduces to the much simpler expression
\begin{equation*}
l_j(x) = \prod_{\substack{1 \leq k \leq N \\ k \neq j}} \frac{x - \xi_k}{\xi_j - \xi_k}, \qquad j=1, \dots, N.
\end{equation*}

The quality of the interpolant $L_n f$ generally depends on the function $f$, as well as the choice of points $X_{n,d}$. We are interested in finding the set $X_{n,d}$ which is optimal in the sense that it gives the lowest possible upper bound on the interpolation error for arbitrary $f \in C^0(D)$. More precisely, let $f_n^* \in \Pi_n^d$ be the best approximation to $f$ in $\Pi_n^d$ in the uniform norm, i.e.
\[
\|f - f_n^*\|_{\infty} = \min_{g \in \Pi_n^d} \|f - g\|_{\infty}.
\] 
Since $f_n^* = L_{n} f_n^*$, we have 
\[
\|f - L_n f\|_{\infty} \, \leq \, \|f - f_n^*\|_{\infty} \, +\, \|L_{n} f_n^* - L_{n} f \|_{\infty} \, \leq \, (1 + \lambda_{n,d}) \; \|f - f_n^*\|_{\infty},
\]
where 
$$
\lambda_{n,d} := \|L_n\|_\infty =  \max_{x \in D} \sum_{j=1}^N |l_j(x)|
$$
is known as the {\em Lebesgue constant} of the set $X_{n,d}$. Note that the above bound holds true for any $f \in C^0(D)$. We then define the optimal points $X_{n,d}^* = \{\xi^*_j\}_{j=1}^N$ as those minimising the Lebesgue constant over all possible choices of interpolation points, and denote the corresponding Lebesgue constant by $\lambda_{n,d}^*$. More precisely, we have 
\begin{equation}\label{eq:minmax}
X_{n,d}^* =  \argmin_{\{\xi_j\}_{j=1}^N} \; \max_{x \in D} \; \sum_{j=1}^N |l_j(x)|, \quad \text{and} \quad \lambda_{n,d}^* = \min_{\{\xi_j\}_{j=1}^N} \; \max_{x \in D} \; \sum_{j=1}^N |l_j(x)|.
\end{equation}

It follows from \eqref{eq:lebpol_d} that the Lebesgue constant $\lambda_{n,d}$, as a function of the interpolation points $X_{n,d}$, is continuous at every point $Y \in D^N$ which is unisolvent \cite{bsv12}. By a suitable redefinition of the Lebesgue constant at the point sets which are not unisolvent, the Lebesgue constant becomes lower--semicontinuous and thus can be shown to have a global minimum on the compact set $D^N$. The global minimum is, however, in general not unique, and several local minima exist in addition to the global minimum.

The aim of this paper is to study the growth of the optimal Lebesgue constant $\lambda_{n,d}^*$, both with respect to the degree $n$ and the dimension $d$, and to examine the structure of the optimal interpolation points $X_{n,d}^*$. To this end, let us briefly review some known results on $X_{n,d}^*$ and $\lambda_{n,d}^*$, as well as the behaviour of $\lambda_{n,d}$ for some well known point sets.

\subsection{Review of univariate interpolation}
In the case of univariate interpolation, the behaviour of the optimal Lebesgue constant $\lambda_{n,1}^*$ has been fully characterised. It has been shown (see e.g. \cite{v90,s06}) that for the canonical interval $D=[-1,1]$, $\lambda_{n,1}^*$ allows the expansion
\begin{equation}\label{eq:leb_opt_1d}
\lambda_{n,1}^* = \frac{2}{\pi} \log (n+1) + \frac{2}{\pi}\left( \gamma + \log \frac{4}{\pi} \right) + \mathcal O\left(\left(\frac{\log \log n}{\log n }\right)^2\right),
\end{equation}
where $\gamma \approx 0.577$ denotes the Euler-Mascheroni constant and $\log$ denotes the natural logarithm. Results for any other bounded interval are readily available through the use of linear transformations.

An explicit formula for the construction of the optimal point set $X_{n,1}^*$, on the other hand, is yet unknown. As already noted earlier in this section, without imposing any additional constraints on the optimisation problem \eqref{eq:minmax}, the minimising point set is not unique. However, in the case of univariate interpolation it has been shown that if two of the interpolation points are fixed to be at the left and right end point of the interval, respectively, a unique minimising set is guaranteed. It is furthermore known that this set is symmetric (see e.g. \cite{h98,k78} and the references therein).

Although the exact location of the optimal point set $X_{n,1}^*$ is unknown, explicit expressions exist for {\em nearly optimal} points. In light of \eqref{eq:leb_opt_1d}, nearly optimal point sets are defined as sets which have a Lebesgue constant that grows logarithmically in $(n+1)$, but with possibly larger constants than $\lambda_{n,1}^*$. 
An example of nearly optimal points are the Gauss-Chebyshev-Lobatto nodes, which on the canonical interval $[-1,1]$ are given by $\xi_1 = 0$ for $N=1$ and otherwise by
\[
\xi_j = -\cos\left(\frac{\pi(j-1)}{N-1}\right), \quad j=1, \dots, N.
\]
Of all known nearly optimal interpolation points, these nodes have the smallest Lebesgue constant \cite{s06,h98}, which is given by
\[
\lambda_{n,1} = \frac{2}{\pi} \log (n+1) + \frac{2}{\pi}\left( \gamma + \log \frac{8}{\pi} -\frac{2}{3} \right) + \mathcal O\left(\frac{1}{\log n}\right).
\]

It is conjectured (see e.g. \cite{cg10}) that nearly optimal point sets are asymptotically equidistributed with respect to the Dubiner metric
\[
d(x_1,x_2) = |\arccos x_2 - \arccos x_1|.
\]
Examples of nearly optimal point sets which have this property include Gauss-Chebyshev, Gauss-Chebyshev-Lobatto and Leja points \cite{cdv05}.

\subsection{Review of multivariate interpolation}
In the case of multivariate interpolation, much less is known about the optimal point sets and the growth of the optimal Lebesgue constant. Theoretical results on the Lebesgue constant are known in only a few special cases, and explicit constructions for nearly optimal point sets for total degree interpolation are, to the best of our knowledge, known only in the case of bivariate interpolation.

The only instance where the growth of the optimal Lebesgue constant is known, seems to be the case of interpolation in the unit ball $D = B^d :=\{ x \in \R^d : \|x\| \leq 1\}$. In this case, $\lambda_{n,d}^*$ has been proven to grow algebraically in $n$ for $d > 1$ \cite{s84}: 
\begin{equation}\label{eq:growth_ball}
c_1^B(d) \, n^{(d-1)/2} \; \leq \; \lambda_{n,d}^* \; \leq \; c_2^B(d) \, n^{(d-1)/2},
\end{equation}
for some constants $c_1^B(d)$ and $c_1^B(d)$ which depend on the dimension $d$, but are independent of the degree $n$. The explicit form of these constants is not known.

Another particular choice of $D$ which is often of interest and easily generalises to the multivariate setting, is that of the unit cube $C^d = [-1,1]^d$. In this case the growth of $\lambda_{n,d}*$ is not known analytically. However, there has been some recent progress on constructing nearly optimal point sets for bivariate interpolation on the square. The so-called Padua points \cite{bcdvx06,cdv05} are an example of an explicitly constructed point set on the square, for which the Lebesgue constant has been proven to grow like the square of the logarithm, 
\[
\lambda_{n,2} \; \leq \; c^P \, (\log (n+1))^2 ,
\]
for some constant $c^P$ independent of $n$. We also mention the Xu points \cite{x96,bdv06}, which are another example of explicitly constructed interpolation points on the square with a Lebesgue constant growing like $(\log (n+1))^2$, albeit using a number of points slightly larger than $N$. Both the Padua and the Xu points are asymptotically equidistributed with respect to the two dimensional Dubiner metric $d(x,y) = \max[|\arccos y_1 - \arccos x_1|,|\arccos y_2 - \arccos x_2|]$, see e.g \cite{cdv05}.

In light of these results, it is conjectured that the Lebesgue constant on the unit cube $C^d$ grows like the $d$th power of the logarithm for $d \geq 1$,
\begin{equation}\label{eq:growth_cube}
\lambda_{n,d}^* \; \leq \; c^C(d) \, (\log (n+1))^d, 
\end{equation}
for some constant $c^C(d)$ independent of $n$. Note that the growth of the Lebesgue constant is in this case much slower than in the case of the unit ball in \eqref{eq:growth_ball}, and that the growth rate in \eqref{eq:growth_cube} is the same as for the Lebesgue constant for tensor product polynomial interpolation (i.e. interpolation in $\bigotimes_{k=1}^d \Pi_n^1$) on the cube $C^d$, see e.g. \cite{bsv12}.

As in the case of univariate interpolation, the optimal points $X_{n,d}^*$ are generally not unique unless one imposes further constraints on the optimisation problem \eqref{eq:minmax}. However, it seems that in the multivariate setting, no geometrical constraint is known which guarantees uniqueness and preserves unisolvence.

\section{Algorithmic considerations}\label{sec:alg}
The numerical computation of the optimal interpolation points $X_{n,d}^*$ requires the solution of the non--linear optimisation problem \eqref{eq:minmax}. The number of variables involved in the optimisation is equal to $d N = d \, \binom{n+d}{d}$, the number of coordinates of the optimal points, and the optimisation process is hence computationally very intensive. 

There are several issues one needs to consider when solving the optimisation problem \eqref{eq:minmax} numerically. 
Firstly, the maximum over $x \in D$ appearing in \eqref{eq:minmax} has to be approximated by the maximum taken over only a finite number of points. Hence, we are looking for finite sets $Y_{n,d} \subset D$, such that 
\begin{equation}\label{eq:leb_adm}
\lambda_{n,d}^* \approx \min_{\{\xi_j\}_{j=1}^N} \; \max_{ x \in Y_{n,d}} \; \sum_{j=1}^N |l_j(x)|.
\end{equation}
In general, the appropriate choice of approximating mesh $Y_{n,d}$ will depend on the degree $n$ and the dimension $d$, and the cardinality of $Y_{n,d}$ should grow with $n$ and $d$. A good choice for $Y_{n,d}$ are {\em admissible meshes}, which were first introduced in \cite{cl08}. For a fixed dimension $d$, an admissible mesh is a sequence $\{Y_{n,d}\}_{n \in \mathbb N}$ of finite subsets of $D$ such that the cardinality of $Y_{n,d}$ grows at most polynomially in $n$ as $n \rightarrow \infty$, and the inequality
\begin{equation}\label{eq:adm_mesh}
\max_{ x \in D} \; |p(x)| \; \leq \; C(Y_{n,d}) \, \max_{x \in Y_{n,d}} \; |p(x)|, \qquad \text{for all} \; p \in \Pi_n^d,
\end{equation}
holds with a constant $C(Y_{n,d})$ bounded from above, independently of $n$. As such, the bound in \eqref{eq:adm_mesh} does not apply directly to the Lebesgue constant $\max_{x \in D} \, \sum_{i=1}^N |l_j(x)|$, since this involves the sum of absolute values of polynomials as opposed to the absolute value of a single polynomial. However, in practice admissible meshes seem to work very well also for approximating the Lebesgue constant $\lambda_{n,d}$ (see section \ref{sec:num}). Other choices of approximating mesh are of course also possible (see e.g. \cite{vhs14}). Examples of admissible meshes are given in section \ref{sec:num}.

Another issue that needs to be addressed, is the efficient evaluation of the Lagrange fundamental polynomials. In order to compute the approximate Lebesgue constant with the admissible mesh $Y_{n,d}$, each of the $N$ Lagrange fundamental polynomials needs to be evaluated at each point in the admissible mesh. One way to do this is to use the explicit expression \eqref{eq:lebpol_d} in terms of Vandermonde determinants. However, for large values of $d$ or $n$, this quickly becomes prohibitively expensive and can in addition suffer from numerical instabilities. An alternative way to assemble the Lagrange fundamental polynomials, for any given set of points $X_{n,d}$, is given in an algorithm by Sauer and Xu \cite[Algorithm 4.1]{sx95}. Starting from any polynomial basis of $\Pi_n^d$, such as the monomial basis, this algorithm assembles the Lagrange fundamental polynomials using only the operations of addition, multiplication by a scalar and point evaluation of polynomials. The algorithm does not require the computation of Vandermonde determinants or the solution of linear systems of equations. Furthermore, proper termination of the algorithm also ensures unisolvence of the points $X_{n,d}$.

Lastly, one has to consider the choice of optimisation algorithm. As already discussed in section \ref{sec:prob}, the Lebesgue constant $\lambda_{n,d}$ is a continuous, and in fact also differentiable, function of the points $X_{n,d}$ everywhere except for the set of non-unisolvent point sets. Due to the presence of these singularities, it is preferable to use algorithms with minimal smoothness assumptions on the objective function. Implementations of various state-of-the-art optimisation algorithms are available through the MATLAB Optimization Toolbox \cite{matopt}. The algorithm most suitable to the problem considered in this paper is the algorithm {\em fminimax}, which is particularly designed to solve min-max optimisation problems and only requires the objective function to be continuous.

\begin{remark}\label{rem:subpol}(Interpolation in subspaces of $\Pi_n^d$) The algorithm described in this section is not specific to interpolation in the total degree spaces $\Pi_n^d$, and can in fact be used to construct optimal point sets for interpolation in general polynomial spaces. Given a set of points $\{\xi_j\}_{j=1}^M$, where $M$ is not equal to the dimension of a total degree polynomial space, the algorithm by Sauer and Xu in \cite{sx95} can be run to compute Lagrange fundamental polynomials $l_1, \dots, l_M$ which span a particular polynomial space, in which the interpolation problem is uniquely solvable, and for which we can then compute optimal interpolation points. 
\end{remark}

The implementation of the optimisation problem \eqref{eq:minmax} has been the subject of several recent studies, such as \cite{vhs14}, \cite{bsv12} and some references therein. However, the focus of these studies has always been the bivariate setting, and higher dimensional problems were not considered. In \cite{vhs14}, the authors present several algorithms to compute minimal (or low) Lebesgue constants based on various approximations in the optimisation problem \eqref{eq:minmax}. All of the algorithms require the solution of linear equations involving generalised Vandermonde matrices. The approach taken by Briani et al in \cite{bsv12} is similar to the approach taken in this paper, although there are substantial algorithmic differences: the Lagrange fundamental polynomials are computed through the use of expression \eqref{eq:lebpol_d}, involving Vandermonde determinants, and a different optimisation algorithm is used ({\em fmincon} from the MATLAB optimisation toolbox, which assumes differentiabilty of the objective function). 

The results obtained in the studies \cite{vhs14} and \cite{bsv12} are quite similar, although the algorithm considered in \cite{vhs14} consistently gives significantly lower Lebesgue constants for large values of $n$. For the range of degrees considered in this paper, the values of the Lebesgue constants reported in \cite{vhs14,bsv12} are always significantly higher or identical to the minimal Lebesgue constants computed in this work.


\section{Numerical investigation}\label{sec:num}
We implemented the optimisation problem \eqref{eq:minmax} as described in section \ref{sec:alg}. The optimisation algorithm used was {\em fminimax} from the MATLAB Optimization toolbox \cite{matopt}, and the objective function supplied to {\em fminimax} was the Lebesgue constant approximated on an admissible mesh, as defined in \eqref{eq:leb_adm}. During each iteration of the optimisation procedure, the Lagrange polynomials were assembled using the algorithm by Sauer and Xu, \cite[Algorithm 4.1]{sx95}.

Examples of admissible meshes in various geometries were considered in \cite{bv12}. 
On the unit cube $C^d$, we use the admissible mesh given by tensor product Gauss-Chebyshev-Lobatto grids: $Y_{n,d}^C = \bigotimes_{k=1}^d \text{GCL}_m$, where $\text{GCL}_m$ denotes the set of $m$ Gauss-Chebyshev-Lobatto points in the unit interval $[-1,1]$. For any $m>n$, inequality \eqref{eq:adm_mesh} is in this case satisfied with constant $C(Y_{n,d}^C) = (\cos(\frac{\pi n}{2 m}))^{-d}$ \cite{bv12}. 

On the unit disk $B^2$, \cite{bv12} provides an admissible mesh defined in (non-standard) polar coordinates, consisting of a radial coordinate $r \in [-1,1]$ and an angular coordinate $\theta \in [0,\pi)$. Compared to using standard polar coordinates, this approach leads to an admissible mesh with lower cardinality for the same constant $C(Y_{n,2}^B)$. The admissible mesh is given by the tensor product of a Gauss-Chebyshev-Lobatto grid in $r$ and a uniform grid in $\theta$: $Y_{n,2}^B = \text{GCL}_m \bigotimes \{\frac{k\pi}{m} : 0 \leq k \leq m-1 \}$. Inequality \eqref{eq:adm_mesh} is again satisfied, with $C(Y_{n,2}^B) = (\cos(\frac{\pi n}{2 m}))^{-2}$. 

Following the tensor product construction, we build an admissible mesh on the three dimensional unit ball $B^3$ in standard spherical coordinates $r \in [0,1]$, $\theta \in [0,\pi]$ and $\phi \in [0,2\pi)$ through the tensor product $Y_{n,3}^B = \text{GCL}_m \bigotimes \{\frac{k\pi}{m} : 0 \leq k \leq m-1 \} \bigotimes \{ \frac{k\pi}{m} : 0 \leq k \leq 2m-1 \}$, where $\text{GCL}_m$ are now $m$ Gauss-Chebyshev-Lobatto points in the interval $[0,1]$. 
For all admissible meshes, we choose the number of points in each dimension, denoted by $m$, growing as $m(n) = 2^{n-1}+1$. 

In an attempt to recover the global minimum of the Lebesgue constant, the initial guess for the optimisation procedure was simply chosen as a random set of points in $D$. The optimisation process was performed over 10 such random initial guesses, and the final Lebesgue constant is taken as the minimum over the 10 computed Lebesgue constants. 

To reduce the computational effort, for each of the 10 starting guesses, the optimisation was performed in an iterative fashion on a sequence of admissible meshes. First, using the random starting guess, the optimisation was performed on a coarse admissible mesh, typically containing about $3^d$ or $5^d$ points. The result of this optimisation was then used as a starting guess for the optimisation on the next admissible mesh. This procedure was continued until refining the admissible mesh resulted in a difference in the Lebesgue constant that was less than $10^{-3}$. Preliminary numerical investigations showed that this iterative approach does not alter the final value found for the optimal Lebesgue constant, i.e. the iterative procedure results in the same minimal Lebesgue constant as an optimisation done directly on the finest admissible mesh.

For the majority of our numerical investigation, we will focus on the case of interpolation in the unit cube $C^d=[-1,1]^d$, and one of our main aims will be to confirm the logarithmic growth rate of the Lebesgue constant as conjectured in \eqref{eq:growth_cube}. We also present some results for interpolation in the unit ball $B^d$.

\subsection{Results in the cube}\label{ssec:num_cube}
We start with the interpolation problem posed on the unit cube $D=[-1,1]^d$. The smallest values of the Lebesgue constant found are given in Table \ref{tbl:leb_cube}. 

\begin{table} [th]
\begin{center}
\renewcommand{\arraystretch}{1.25}
\begin{tabular}{ |c|| c c c c c c c c c c|} \hline
\backslashbox{$d$}{$n$}& 1& 2& 3& 4& 5& 6& 7& 8& 9& 10 \\ \hline \hline
1&     1& 1.25& 1.42& 1.56& 1.67& 1.77& 1.85& 1.93& 1.99& 2.05 \\
2&  1.89& 2.38& 2.73&  3.12&  3.51& 3.86& 4.18& 4.44& 4.71& 4.96 \\
3&  2.00& 2.95& 4.05&  5.09&  6.40& & & & & \\ 
4&  2.32& 3.67& 5.40&  7.85& & & & & & \\
5&  2.45& 4.20& &  & & & & & & \\
6& 2.60& 5.32& & & & & & & & \\
7& 2.74& 6.55& & & & & & & & \\ 
8& 3.10& 7.74& & & & & & & & \\ 
9& 3.25& & & & & & & & & \\
10& 3.56& & & & & & & & & \\ \hline
\end{tabular}
\end{center}
\caption{Optimal Lebesgue constant $\lambda_{n,d}^*$ for various values of $n$ and $d$, for $D =[-1,1]^d$.}
\label{tbl:leb_cube}
\end{table}

\begin{figure}[p]
\centering
\includegraphics[width=0.45\textwidth]{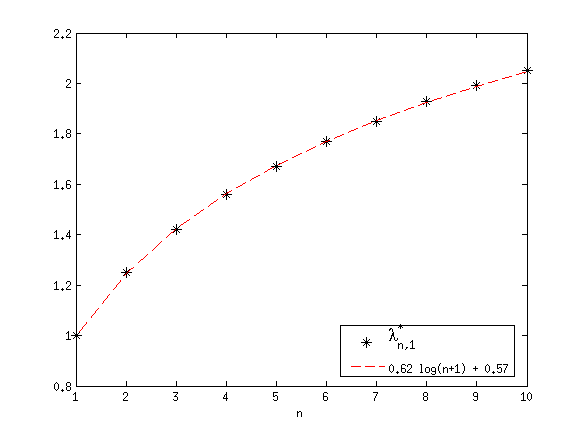}\includegraphics[width=0.45\textwidth]{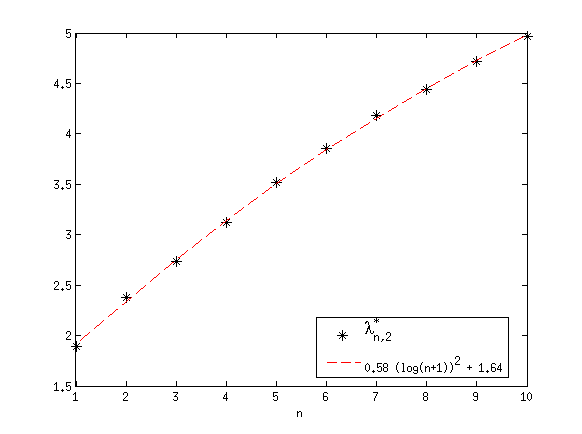}\\
\includegraphics[width=0.45\textwidth]{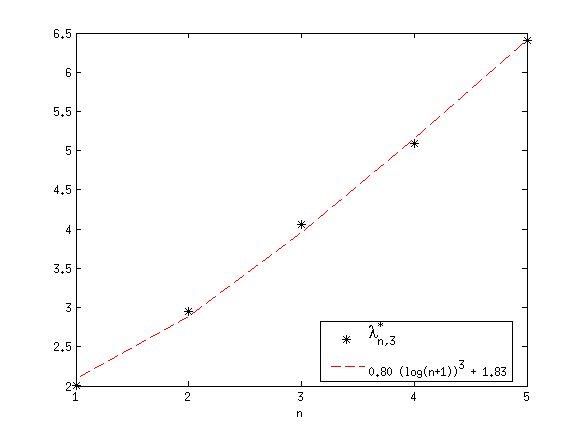}\includegraphics[width=0.45\textwidth]{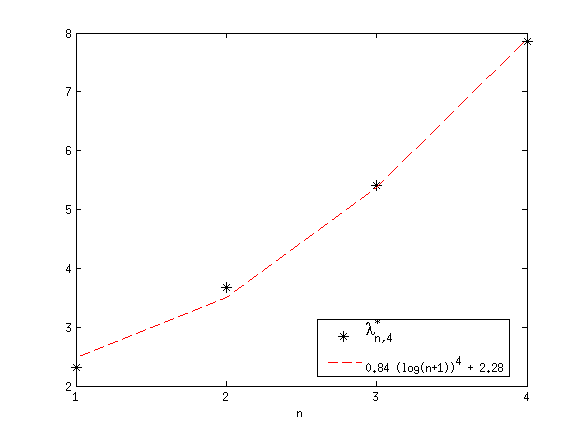}
\caption{Optimal Lebesgue constant $\lambda_{n,d}^*$ as a function of $n$ for $d=1$ (top left), $d=2$ (top right), $d=3$ (bottom left) and $d=4$ (bottom right), for $D =[-1,1]^d$.}
\label{fig:leb_cube_d}
\end{figure}

Let us examine these results in more detail. In particular, we are interested in quantifying the growth of the Lebesgue constants reported in Table \ref{tbl:leb_cube} with respect to both dimension $d$ and degree $n$. In Figure \ref{fig:leb_cube_d}, we have plotted the Lebesgue constant $\lambda_{n,d}^*$ as a function of $n$, for $d=1,2,3,4$. In each of the graphs, we have added a best fit of the growth, which in each case turned out to be of the form
\[
\lambda_{n,d}^* \approx c_1^C(d) (\log(n+1))^d \; + \; c_2^C(d).
\]
This confirms the growth conjectured in \eqref{eq:growth_cube}. The coefficients $c_1^C(d)$ and $c_2^C(d)$ were, for each $d$ individually, determined by a least squares approach.

\begin{figure}[p]
\centering
\includegraphics[width=0.45\textwidth]{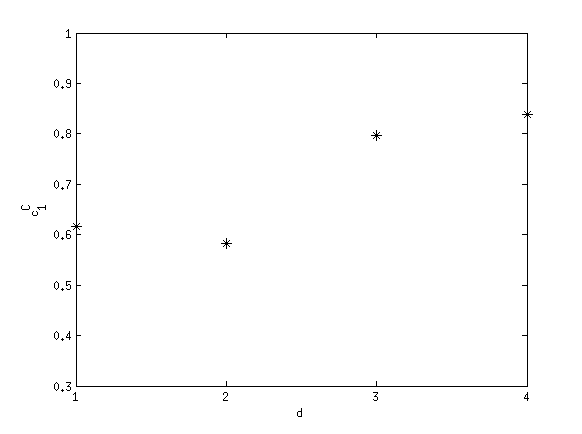} \includegraphics[width=0.45\textwidth]{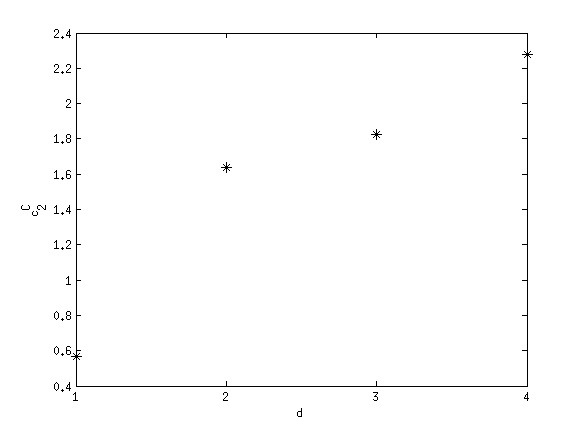}
\caption{Coefficients $c_1^C(d)$ (left) and $c_2^C(d)$ (right) as a function of $d$, for $D =[-1,1]^d$.}
\label{fig:c1_c2}
\end{figure}

It is clear from Figure \ref{fig:leb_cube_d} that the constants $c_1^C(d)$ and $c_2^C(d)$ grow with dimension $d$. Using the values depicted in Figure \ref{fig:leb_cube_d}, we have the estimates of $c_1^C(d)$ and $c_2^C(d)$ shown in Figure \ref{fig:c1_c2}. It appears that for $d \geq 2$, the coefficient $c_1^C(d)$ grows slower than linearly in $d$, whereas $c_2^C(d)$ seems to grow linearly or faster. Further investigations are required to characterise the behaviour of these coefficients more precisely.

\begin{figure}[p]
\centering
\includegraphics[width=0.45\textwidth]{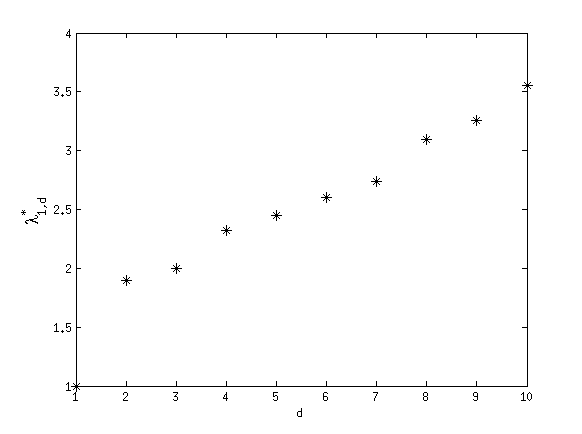} \includegraphics[width=0.45\textwidth]{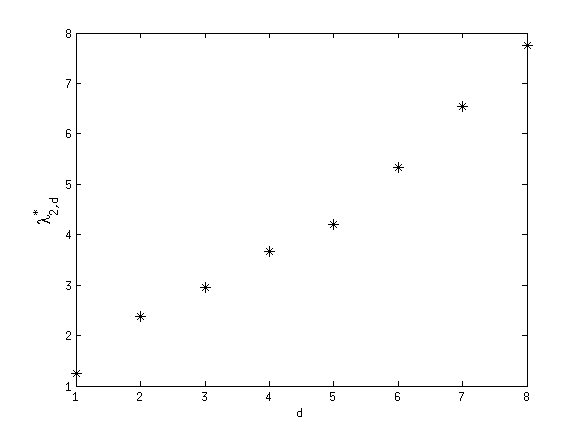}
\caption{Optimal Lebesgue constant $\lambda_{n,d}^*$ as a function of $d$ for $n=1$ (left) and $n=2$ (right), for $D =[-1,1]^d$.}
\label{fig:leb_cube_n}
\end{figure}

In Figure \ref{fig:leb_cube_n}, we have plotted the optimal Lebesgue constant $\lambda_{n,d}^*$ as a function of the dimension $d$, for $n=1,2$. It appears that $\lambda_{1,d}^*$ grows linearly in $d$, whereas $\lambda_{2,d}^*$ grows somewhat faster. This is again in accordance with the suggested growth rate $\lambda_{n,d}^* \leq c^C (\log(n+1))^d$ in \eqref{eq:growth_cube}, which for $\log(n+1)>1$, i.e. $n \geq 2$, predicts an exponential growth in $d$.

\begin{figure}[p]
\centering
\includegraphics[width=0.5\textwidth]{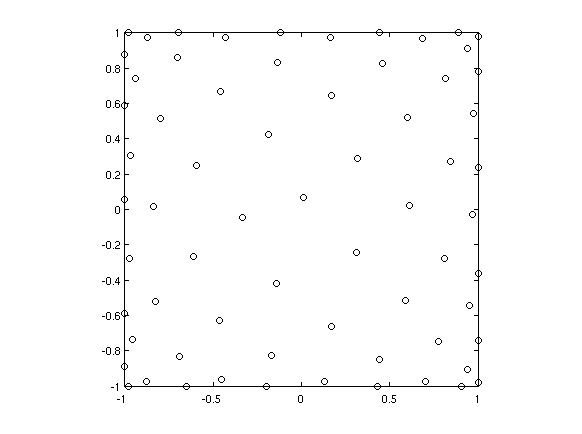}\includegraphics[width=0.5\textwidth]{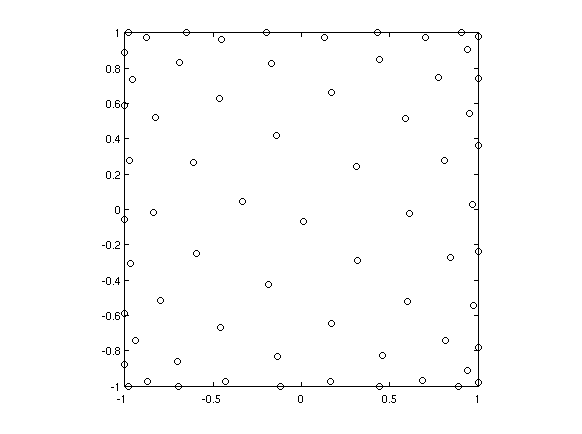}\\
\caption{$N=66$ optimal interpolation points for $d=2$ and $n=10$.}
\label{fig:opt_2d}
\end{figure}

Let us now move on to studying the structure of the optimal interpolation points $X_{n,d}^*$. Two possible choices of the set $X_{10,2}^*$ are shown in Figure \ref{fig:opt_2d}. We see that although the points are spaced quite evenly throughout the square, the density of points is closer to the boundary. As mentioned in section \ref{sec:prob}, the minimising sets $X_{n,2}^*$ are in general not unique. In our numerical experiments we found, however, that all the minimising sets can be obtained from one minimising set through reflections in the $x_1$ axis, $x_2$ axis and/or the diagonals of the square. The Lebesgue constant is unaffected by these transformations, and which of these configurations is found by the optimisation algorithm depends on the initial guess. The two sets in Figure \ref{fig:opt_2d} are related through a reflection in the horizontal $x_1$ axis. 

\begin{figure}[p]
\centering
\hspace*{-5ex}\includegraphics[width=0.5\textwidth]{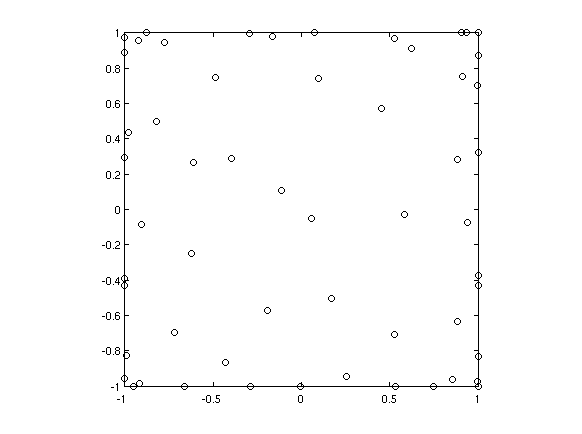} \includegraphics[width=0.5\textwidth]{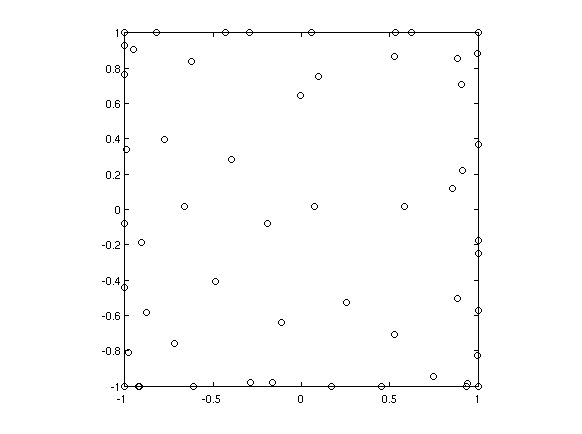} \includegraphics[width=0.5\textwidth]{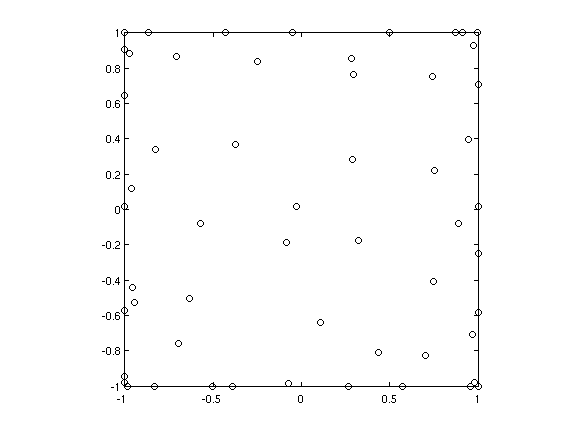}
\caption{Two dimensional projections of $N=56$ optimal interpolation points for $d=3$ and $n=5$. Top left: $x_1$-$x_2$. Top right: $x_1$-$x_3$. Bottom: $x_2$-$x_3$.}
\label{fig:opt_3d}
\end{figure}

\begin{figure}[p]
\centering
\hspace*{-5ex}\includegraphics[width=0.5\textwidth]{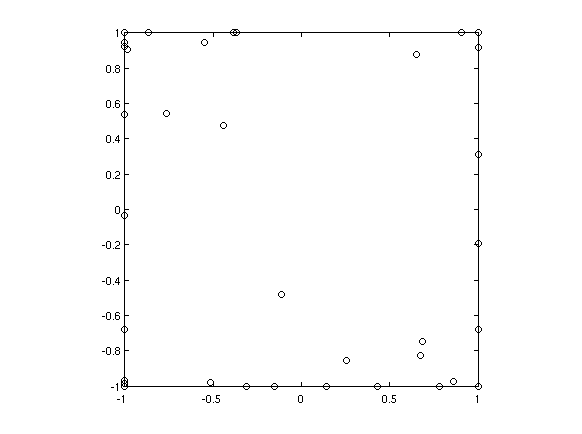} \includegraphics[width=0.5\textwidth]{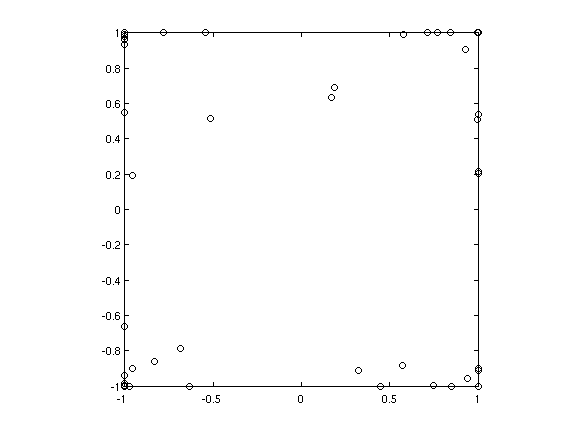} 
\caption{Two dimensional projections of $N=45$ optimal interpolation points for $d=8$ and $n=2$. Left: $x_1$-$x_2$. Right: $x_7$-$x_8$.}
\label{fig:opt_8d}
\end{figure}

The same trends are evident in optimal interpolation points in higher dimensions. Reflections in one or more of the coordinate axes or diagonals do not change the value of the Lebesgue constant, and the density of points is higher towards the boundary of the domain. For small values of $n$, the points seem to be located exclusively around the boundary, and the centre of the cube gets filled in as $n$ increases. The two dimensional projections of an optimal point set $X_{5,3}^*$ are shown in Figure \ref{fig:opt_3d}, and a selection of two dimensional projections of an optimal point set $X_{2,8}^*$ are shown in Figure \ref{fig:opt_8d}.

Lastly, let us compare our results found in the case $d=2$ to those found by Briani et al in \cite{bsv12} and by van Barel et al in \cite{vhs14}. The specific values of the Lebesgue constants are shown in Table \ref{tbl:comp_cube}, and confirm a significant improvement. Although only a graph of the Lebesgue constants is given by van Barel et al in \cite{vhs14}, without mention of the exact values, it is clear from \cite[Figure 6.1]{vhs14} that their results are almost identical to those found in \cite{bsv12}, and the Lebesgue constants computed in this paper are hence significantly lower also than the ones computed in \cite{vhs14}.

\begin{table} [h]	
\begin{center}
\renewcommand{\arraystretch}{1.25}
\begin{tabular}{ |c|| c c c c c c c c c c|} \hline
$n$ & 1& 2& 3& 4& 5& 6& 7& 8& 9& 10 \\ \hline \hline
Briani et al & 2.00& 2.39& 2.73& 3.24& 3.59& 4.00& 4.34& 4.90& 5.18& 5.32 \\ 
Our results  & 1.89& 2.38& 2.73& 3.12& 3.51& 3.86& 4.18& 4.44& 4.71& 4.96\\\hline
\end{tabular}
\end{center}
\caption{Optimal Lebesgue constants $\lambda_{n,2}^*$, for $D=[-1,1]^2$.}
\label{tbl:comp_cube}
\end{table}

\begin{figure}[p]
\centering
\hspace{-5ex}\includegraphics[width=0.5\textwidth]{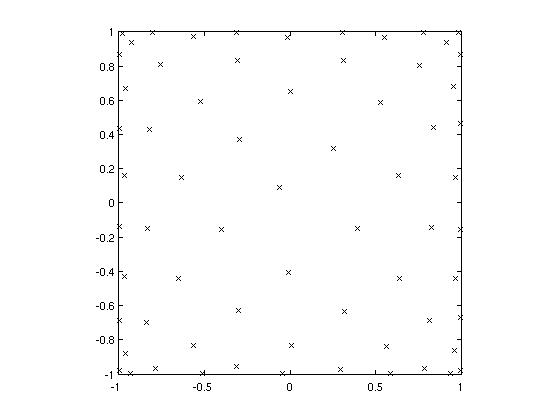} \includegraphics[width=0.5\textwidth]{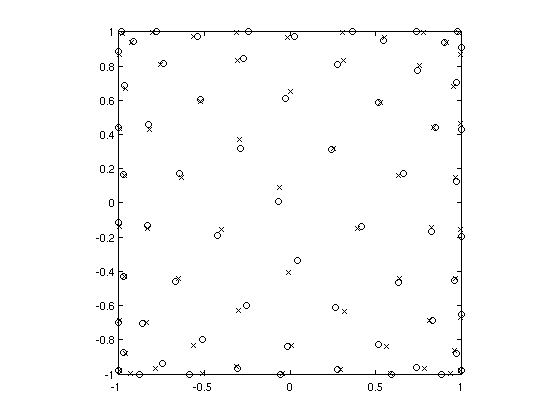}
\caption{Comparison of optimal points $X_{10,2}^*$. Left: points computed by Briani et al in \cite{bsv12}. Right: comparison of the points found in \cite{bsv12} and the points found in this paper.}
\label{fig:d2n10_comp}
\end{figure}

The optimal point sets $X_{10,2}^*$ as found in this paper and by Briani et al in \cite{bsv12} are shown in Figure \ref{fig:d2n10_comp}. The left plot shows only the points computed in \cite{bsv12}, whereas the right plot shows both sets of points. For a fair comparison, we have chosen the configuration of our optimal points that best matches the Briani et al set. Although the two point sets appear very similar, significant differences are visible, especially in the centre of the square.

\subsection{Results in the ball}\label{ssec:num_ball}
We now repeat some of the experiments in the previous section in the unit ball $D = \{ x \in \R^d : \|x\| \leq 1\}$. The smallest values of the Lebesgue constant found are reported in Table \ref{tbl:leb_ball}.

\begin{table} [th]
\begin{center}
\renewcommand{\arraystretch}{1.25}
\begin{tabular}{ |c|| c c c c c c c c c c|} \hline 
\backslashbox{$d$}{$n$}& 1& 2& 3& 4& 5& 6& 7& 8& 9& 10 \\ \hline \hline
2&  1.67& 1.99& 2.47&  2.95&  3.39& 3.85& 4.30& 4.84& 5.20& 5.69\\
3&  2.00& 3.06& 3.56& 4.78& 5.92& & & & & \\ \hline
\end{tabular}
\end{center}
\caption{Optimal Lebesgue constant $\lambda_{n,d}^*$ for various values of $n$ and $d$, for $D = \{ x \in \R^d : \|x\| \leq 1\}$.}
\label{tbl:leb_ball}
\end{table}

\begin{figure}[p]
\centering
\includegraphics[width=0.45\textwidth]{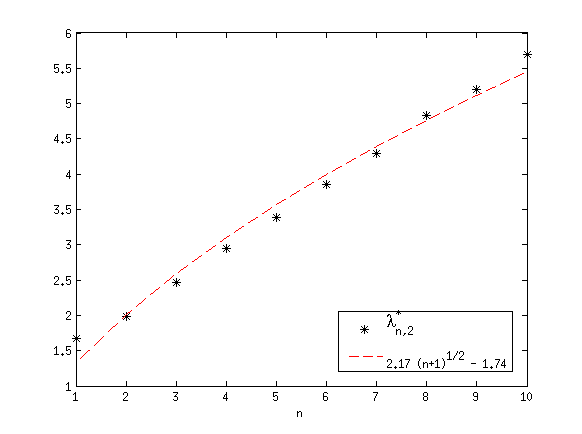} \includegraphics[width=0.45\textwidth]{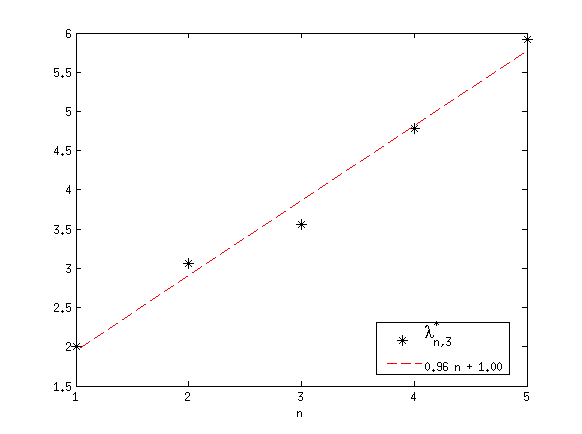}
\caption{Optimal Lebesgue constant $\lambda_{n,d}^*$ as a function of $n$ for $d=2$ (left) and $d=3$ (right), for $D = \{ x \in \R^d : \|x\| \leq 1\}$.}
\label{fig:leb_n_ball}
\end{figure}

The growth of the Lebesgue constants $\lambda_{n,d}^*$ with respect to $n$, for $d=2,3$, is shown in Figure \ref{fig:leb_n_ball}. We have again added a least squares fit to the graphs. For $d=2$, the Lebesgue constant seems to grow like the square root of $(n+1)$, whereas for $d=3$ it grows linearly in $n$. These results are in agreement with the theoretical result in \eqref{eq:growth_ball}.

\begin{figure}[p]
\centering
\hspace{-5ex}\includegraphics[width=0.5\textwidth]{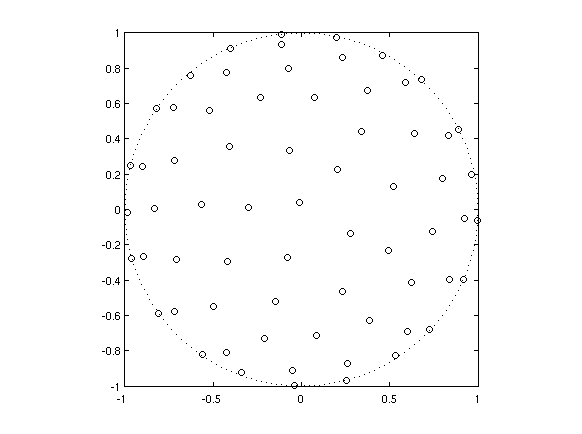} \includegraphics[width=0.5\textwidth]{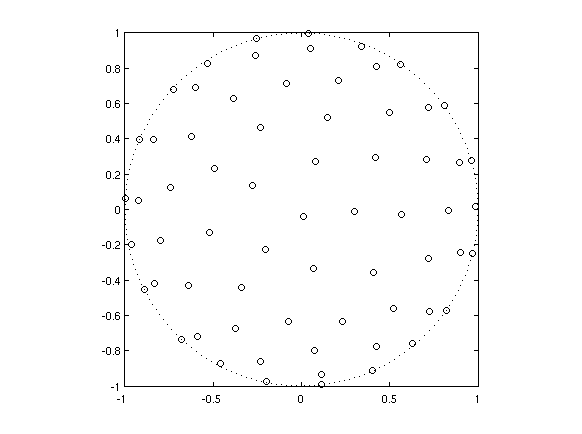}
\caption{$N=66$ optimal interpolation points for $d=2$ and $n=10$.}
\label{fig:d2n10_disk}
\end{figure}

\begin{figure}[p]
\centering
\hspace*{-5ex}\includegraphics[width=0.5\textwidth]{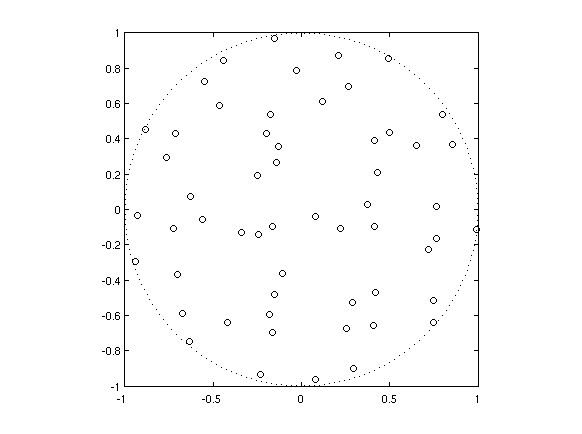} \includegraphics[width=0.5\textwidth]{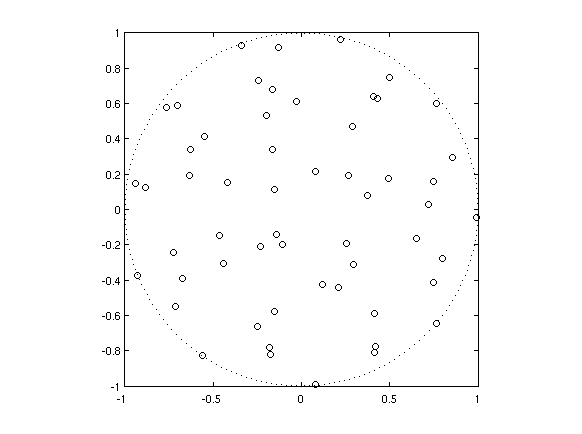} \\ \hspace*{-11ex} \includegraphics[width=0.5\textwidth]{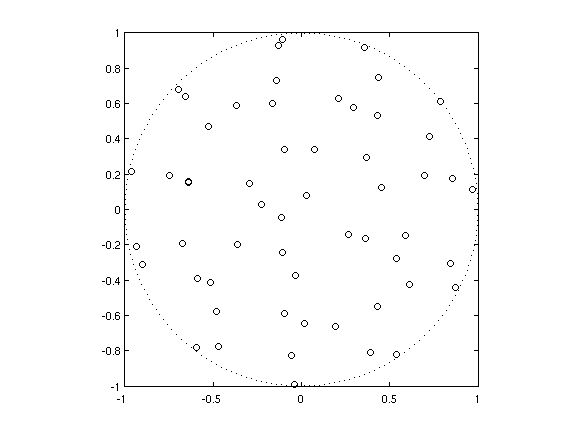} \includegraphics[width=0.45\textwidth]{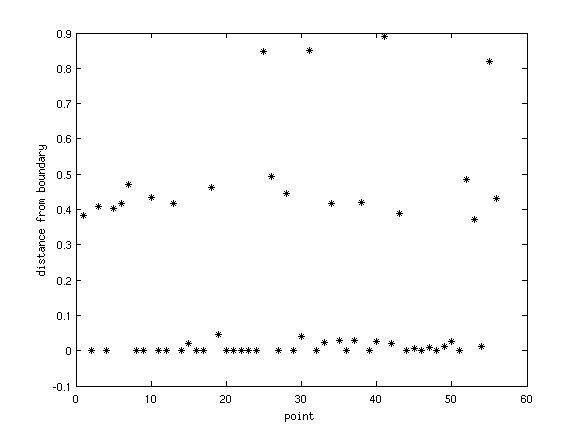}
\caption{Two dimensional projections of $N=56$ optimal interpolation points for $d=3$ and $n=5$, together with the Euclidean distance from the boundary $\|x\|=1$ for each of the $56$ points. Top left: $x_1$-$x_2$. Top right: $x_1$-$x_3$. Bottom left: $x_2$-$x_3$. Bottom right: Distance from boundary.}
\label{fig:opt_3d_ball}
\end{figure}

Two configurations of optimal point sets $X_{10,2}^*$ are shown in Figure \ref{fig:d2n10_disk}. As already established for the square in the previous section, the density of points is again higher towards the boundary of the disk, and the optimal points are invariant under reflections through the coordinate axes. Additionally, the Lebesgue constant appears to be invariant under any rotation. The two point sets in Figure \ref{fig:d2n10_disk} are related through a reflection in the vertical $x_2$ axis. The two dimensional projections of an optimal point set $X_{5,3}^*$ are shown in Figure \ref{fig:opt_3d_ball}. To confirm that we again see a higher density of points near the boundary, we have have also plotted the distance from the boundary $\|x\|=1$. 

A comparison of the minimal Lebesgue constants found in this paper and by Briani et al in \cite{bsv12} is given in Table \ref{tbl:comp_ball}. A clear improvement is visible. Specific values of the Lebesgue constants are again not given in \cite{vhs14}, but it is clear from \cite[Figure 6.1]{vhs14} that the Lebesgue constants found in this paper are significantly lower also than those found in \cite{vhs14}.

\begin{table} [th]	
\begin{center}
\renewcommand{\arraystretch}{1.25}
\begin{tabular}{ |c| c c c c c c c c c c|} \hline
$n$ & 1& 2& 3& 4& 5& 6& 7& 8& 9& 10 \\ \hline
Briani et al & 1.67& 1.99& 2.47& 2.97& 3.50& 4.30& 5.08& 5.43& 6.73& 7.62 \\ 
Our results  & 1.67& 1.99& 2.47& 2.95& 3.39& 3.85& 4.30& 4.84& 5.20& 5.69 \\\hline
\end{tabular}
\end{center}
\caption{Optimal Lebesgue constants $\lambda_{n,2}^*$, for $D=\{ x \in \R^2 : \|x\| \leq 1\}$.}
\label{tbl:comp_ball}
\end{table}

\begin{figure}[p]
\centering
\hspace{-5ex}\includegraphics[width=0.5\textwidth]{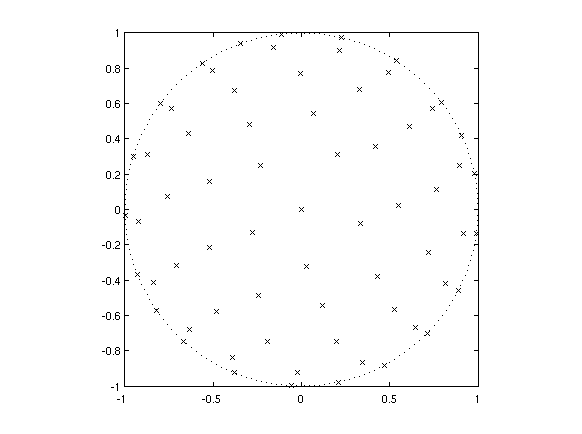} \includegraphics[width=0.5\textwidth]{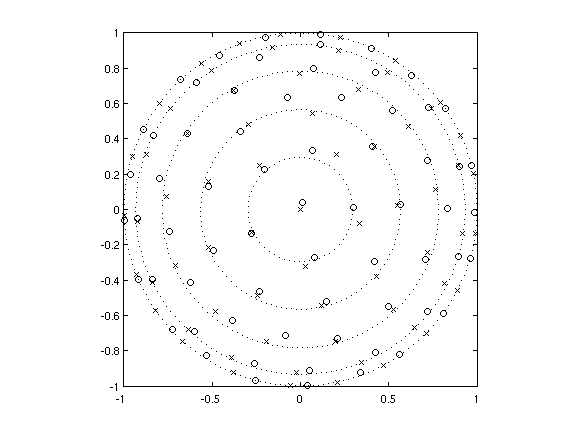}
\caption{Comparison of optimal points $X_{10,2}^*$. Left: points computed in this paper. Right: points computed by Briani et al in \cite{bsv12}.}
\label{fig:d2n10_comp_disk}
\end{figure}

A comparison of the optimal points sets $X_{10,2}^*$ as found in this paper and by Briani et al in \cite{bsv12} is shown in Figure \ref{fig:d2n10_comp_disk}. The left plot shows only the points computed in \cite{bsv12}, whereas the right plot shows both sets of points. Differences between the two point sets are clearly visible, although the structure of the two sets appears similar. It was suggested in \cite{bsv12} that the optimal points (approximately) lie on concentric circles with radii distributed as positive Gauss-Legendre-Lobatto points, so we have added these for comparison. Compared to the points found in \cite{bsv12}, the points found in this paper are further away from this pattern. Note also that a similar pattern was observed in the bottom left plot in Figure \ref{fig:opt_3d_ball}, where the points seemed to lie approximately on concentric spheres.

\section{Conclusions and further work}\label{sec:conc}
The determination of good interpolation points for multivariate Lagrange interpolation in total degree polynomial spaces remains an open question. In this work, we have addressed the question of characterising the optimal choice of interpolation points, which result from a minimisation of the Lebesgue constant. We provided an algorithm to numerically compute the optimal interpolation points, which does not require the computation of Vandermonde determinants or the solution of linear systems of equations, and we provided optimal sets of interpolation points, together with their Lebesgue constant, in various multivariate settings. The Lebesgue constants reported in this work are, to the best of our knowledge, the lowest known values to date.

In future work, it would be interesting to apply the algorithm in this work to even higher dimensions and higher degree polynomial interpolation. To make the algorithm more computationally feasible in these situations, a major advantage would be the use of admissible meshes which are not of tensor-product type. However, it seems that no examples of such admissible meshes are known. A good initial guess for the interpolation points would also be of great benefit, although care has to be taken in the choice of initial configuration due to the many local minima of the Lebesgue constant.

\section*{Acknowledgements} The authors are supported by the US Department of Energy Advanced Simulation Computing Research (ASCR) program under grants DE-SC0010678 and DE-SC0009324.

\bibliographystyle{plain}
\bibliography{intbib}{}

\end{document}